\documentclass[preprint,12pt]{elsarticle}

\usepackage{amsfonts}
\usepackage{graphicx}
\usepackage{epsfig}
\usepackage{geometry}
\usepackage{tikz}
\usetikzlibrary{trees,positioning,shapes}
\usepackage{subcaption}
\usepackage{amsmath}
\usepackage{amsthm}

\usepackage{xcolor, colortbl}
\usepackage{rotating}

\DeclareSymbolFont{largesymbolsA}{U}{txexa}{m}{n} %
\DeclareMathSymbol{\bigtimes}{\mathop}{largesymbolsA}{16} %

\newtheorem{theorem}{Theorem}

\newtheorem{remark}{Remark}

\newcommand{\Prod}[2]{(#1, #2)_{L^2}}

\newcommand{\Vect}[1]{\mathbf{#1}}

\newcommand{\SGM}{\boldsymbol{\sigma}}

\newcommand{\bfs}[1]{{\boldsymbol #1}}

\begin{document}

\begin{frontmatter}

\title{Fast isogeometric solvers for hyperbolic wave propagation problems}

\author{M. \L{}o\'{s},$^{{\textrm{(1)}}}$, P. Behnoudfar$^{{\textrm{(2)}}}$,  M. Paszy\'{n}ski $^{{\textrm{(1)}}}$, V. M. Calo$^{{\textrm{(2,3,4)}}}$}

\address{$^{\textrm{(1)}}$ Department of Computer Science, Faculty of Computer Science, Electronics and Telecommunications, AGH University of Science and Technology,
Krakow, Poland \\
e-mail: paszynsk@agh.edu.pl \\
e-mail: marcin.los.91@gmail.com}

\address{$^{\textrm{(2)}}$ 
Applied Geology Department, \\ Western Australian School of Mines, \\
Faculty of Science and Engineering, Curtin University, Perth, WA, Australia,\\
e-mail: victor.calo@curtin.edu.au}

\address{$^{\textrm{(3)}}$  Mineral Resources, Commonwealth Scientific and Industrial Research Organisation (CSIRO), \\ Kensington, WA, Australia 6152}

\address{$^{\textrm{(4)}}$ 
Curtin Institute for Computation, Curtin University, Perth, WA, Australia 6845}

\begin{abstract}

We use the alternating direction method to simulate implicit dynamics.  Our spatial discretization uses isogeometric analysis. Namely, we simulate a (hyperbolic) wave propagation problem in which we use tensor-product B-splines in space and an implicit time marching method to fully discretize the problem.  We approximate our discrete operator as a Kronecker product of one-dimensional mass and stiffness matrices. As a result of this algebraic transformation, we can factorize the resulting system of equations in linear (i.e., $O(N)$) time at each step of the implicit method. We demonstrate the performance of our method in the model P-wave propagation problem. We then extend it to simulate the linear elasticity problem once we decouple the vector problem using alternating triangular methods.  We proof theoretically and experimentally the unconditional stability of both methods.

\end{abstract}
	

\end{frontmatter}

\section{Introduction}

The alternating directions method (ADS) introduced in~\cite{ ADS1, ADS2, ADS3, ADS4} to deal with finite difference simulations for time-dependent problems.  The method currently solves a broad class of problems~\cite{ Minev1, Minev2}.

Isogeometric analysis (IGA)~\cite{ IGA}, uses B-splines or NURBS~\cite{ NURBS} basis functions in finite element simulations. IGA has multiple applications in time-dependent simulations, including phase field models~\cite{ Dede2010, Dede2011}, phase-separation simulations with application to cancer growth simulations~\cite{ Gomez2008, Gomez2010}, wind turbine aerodynamics~\cite{ Hsu:2011}, incompressible hyper-elasticity~\cite{ Duddu:2011}, turbulent flow simulations~\cite{ Chang:2011}, transport of drugs in cardiovascular applications~\cite{ Hossain:2011} as well as the blood flow simulations and drug transport in arteries simulations~\cite{ Bazilevs2006, Bazilevs2007, Calo:2008}.

Recently, Gao et al.~\cite{ CG1, CG2, Gao2014} applied the direction splitting method to the rapid solution of explicit dynamics using isogeometric analysis on tensor-product grids.  These direction splitting schemes deliver a fast inversion method for the spatial discretization by grouping one-dimensional B-splines along particular spatial axes. For similar methods for fast simulations of explicit dynamics refer to~\cite{ ED1, ED2, ED3, ED4, ED5, advances}.

In this paper, we extend this methodology to hyperbolic scalar problems by collecting different terms as a sequence of multi-banded inversions. Then, we extend these ideas to hyperbolic vector problems, where the model problem is isotropic linear elasticity. First, its corresponding differential operator is decoupled (for more details, see~\cite{ samarskii1964economical}), and then our idea is employed. Finally, we prove the unconditional stability of the schemes as well as their order of convergence.

The structure of the paper is the following. In section 2, we start from the description of the direction splitting for the P-wave equation. Next, in Section 3, we show the stability analysis for the P-wave problem. Section 4 presents the numerical results for the three-dimensional P-wave propagation problem. In Section 5, we extend our method to elastic wave propagation, stability analysis in Section 6, and provide numerical evidence in Section 7. In Section 8, we analyze the order of the schemes. We describe our conclusions in Section 9.

\section{Direction splitting for scalar P-wave equation}

We describe the methodology by directly applying it to a model problem. We first solve the scalar P-wave equation problem given by
\begin{equation}
\ddot{u} - \Delta u = f,
\end{equation}
where the over dot represents a time derivative, and $\Delta$ is the  Laplacian operator. We discretize time as follows
\begin{equation}
\ddot{u}_{n+1} - \Delta u_{n+1} = f_{n+1}
\end{equation}
and use a Newmark expansion from time step $n$ to $n+1$~\cite{ Newmark}
\begin{equation}
  u_{n+1} = u_n + \tau\,\dot{u_n} + \frac{1}{4}\tau^2 \,\ddot{u}_{n+1}
\end{equation}
so that
\begin{equation}
\ddot{u}_{n+1} - \frac{1}{4}\tau^2\,\Delta \ddot{u}_{n+1} = \Delta u_n + \tau\, \Delta \dot{u}_n + f_{n+1}
\end{equation}
We treat~$u_n$, $\dot{u}_n$ and~$\ddot{u}_n$ as three independent variables. Thus, we can update~$\dot{u}_n$ according to
\begin{equation}
  \dot{u}_{n+1} = \dot{u}_n + \frac{1}{2} \tau\, \ddot{u}_{n+1}
\end{equation}
As for the~$u_n$, we use a backward Taylor expansion to obtain
\begin{equation}
  u_n = u_{n+1} - \tau\, \dot{u}_{n+1} + \frac{1}{2}\tau^2\,\ddot{u}_{n+1}
\end{equation}
and so
\begin{equation}
  u_{n + 1} = u_n + \tau\, \dot{u}_{n+1} - \frac{1}{2}\tau^2\,\ddot{u}_{n+1}
\end{equation}
The full scheme is thus the following:
\begin{equation}
\label{eq:scheme}
\left\{
\begin{aligned}
\ddot{u}_{n+1} - \frac{1}{4}\tau^2\,\Delta \ddot{u}_{n+1} &=
\Delta u_n + \tau\, \Delta \dot{u}_n + f_{n+1} \\
\dot{u}_{n+1} &= \dot{u}_n +  \frac{1}{2}\tau\, \ddot{u}_{n+1} \\
u_{n + 1} &= u_n + \tau\, \dot{u}_{n+1} - \frac{1}{2}\tau^2\,\ddot{u}_{n+1}
\end{aligned}
\right.
\end{equation}
We can compute~$u_{n+1}$ and~$\dot{u}_{n+1}$ given~$\ddot{u}_{n+1}$.

For the first equation, we test with function $w$. Thus, the full scheme becomes
\begin{equation}
(w,\ddot{u}_{n+1})+ \frac{1}{4} \tau^2 (\nabla w, \nabla \ddot{u}_{n+1}) = (w,\Delta u_n + \tau\, \Delta \dot{u}_n + f_{n+1}).
\end{equation}

We discretize
\begin{equation}
w = \sum_{ab} {n_am_b C_{ab}} \quad 
\ddot{u}_{n+1} = \sum_{cd} {n_cm_d D_{cd}} 
\end{equation}
where $n_a m_b$ and $n_c m_d$ denotes the tensor-product of one-dimensional B-spline, which form a two-dimensional basis function, and $C_{ab}$ and $D_{cd}$ denotes the coefficients associated with the degrees of freedom. The left-hand side of the equation is
\begin{eqnarray}
(n_am_b C_{ab},n_cm_dD_{cd})+\nonumber 
\frac{1}{4} \tau^2(n_a'm_b C_{ab}, n'_cm_d D_{cd}) +
\frac{1}{4} \tau^2(n_am_b' C_{ab}, n_cm'_d D_{cd}).
\end{eqnarray}
Assuming that the geometry of the domain is simple, we can express the mapping as a separable function. Thus, we can now split the left-hand side of the system as follows
\begin{eqnarray}
(n_am_b,n_cm_d)C_{ab}D_{cd}+\nonumber 
(n_a', n'_c)_x\frac{1}{4} \tau^2(m_b,m_d)_yC_{ab}D_{cd} +
(m_b',m'_d)_x\frac{1}{4} \tau^2(m_b',m'_d)_yC_{ab}D_{cd}.
\end{eqnarray}
We define the following one-dimensional mass and stiffness matrices
\begin{eqnarray}
(n_a,n_c)_x=M_x, \nonumber \\
(n_b,n_d)_y=M_y, \nonumber \\
(n'_a,n'_c)_x=K_x, \nonumber \\
(n'_b,n'_d)_y=K_y, 
\end{eqnarray}
and rewrite the entire system as
	\begin{equation}
	\begin{aligned}
\left( M_x \otimes M_y +\frac{1}{4} \tau^2 (M_x \otimes K_y + K_x \otimes M_y)\right) \ddot{U}^{n+1} &= 
-(M_x \otimes K_y + K_x \otimes M_y)U^{n}\\
&-\tau (M_x \otimes K_y + K_x \otimes M_y)\dot{U}^{n}.
+F^{n+1}
\end{aligned}
\end{equation}
We can now approximate the system as
	\begin{equation}
\begin{aligned}
( M_x + \frac{1}{4} \tau^2 K_x) \otimes ( M_y + \frac{1}{4}\tau^2 K_y)  &= 
M_x\otimes M_y \\&+ \frac{1}{4} \tau^2 M_x\otimes K_y + \frac{1}{4} \tau^2 K_x \otimes M_y + \frac{1}{16} \tau^4 K_x \otimes K_y \\&\approx \nonumber 
M_x \otimes M_y +\frac{1}{4} \tau^2 (M_x \otimes K_y + K_x \otimes M_y).
\end{aligned}
\end{equation}
Dropping the red term results in the following
\begin{eqnarray}
\label{eq:pwave2}
\begin{aligned}
( M_x + \frac{\tau^2}{4}  K_x) \otimes ( M_y + \frac{\tau^2}{4}  K_y)  \ddot{U}_{n+1} &= 
-\left(M_x \otimes K_y + K_x \otimes M_y\right)(U^n+\tau\dot{U}^{n})
+F_{n+1} \\
\dot{U}_{n+1} &= \dot{U}_n +\frac{1}{2} \tau\, \ddot{U}_{n+1} \\
U_{n + 1} &= U_n + \tau\, \dot{U}_{n+1} - \frac{1}{2}\tau^2\,\ddot{U}_{n+1}.
\end{aligned}
\end{eqnarray}

\section{Spectral analysis of splitting for wave-propagation problem} \label{sec:ea}
In this section, we analyze the stability of the splitting scheme to show it is unconditionally stable. The analysis follows closely the approach introduced in~\cite{ behnoudfar2018variationally, chung1993time, deng2019high}. Throughout this section, we set $F=0$.

\subsection{Stability of the splitting scheme}\label{stb}
We consider the spectral decomposition of each of the directional matrices $K_\xi$ with respect to its directional $M_\xi$ (see for example~\cite{ horn1990matrix}) and obtain
\begin{equation} \label{eq:sd0}
K_\xi = M_\xi P_\xi D_\xi P_\xi^{-1},
\end{equation}
where $D_\xi$ is a diagonal matrix containing the eigenvalues of the generalized eigenvalue problem
\begin{equation} \label{eq:eigKM0}
K_\xi v_\xi = \lambda_\xi M_\xi v_\xi
\end{equation}
and the columns of $P_\xi$ are the eigenvectors of the generalized problem. Herein, $\xi = x,y,z$ specifies each of the coordinate directions. We state the analysis for 2D splitting and calculate the required terms given by (for details see~\cite{ behnoudfar2018variationally, behnoudfar2019higher})
\begin{equation}\label{eq:multi}
\begin{aligned}
\tilde G^{-1} &= P_x E_x P_x^{-1} M_x^{-1} \otimes P_y E_y P_y^{-1} M_y^{-1},\\
\tilde G^{-1} M &= \left( P_x \otimes P_y \right) \cdot \left( E_x \otimes E_y \right) \cdot \left( P_x^{-1} \otimes P_y^{-1} \right), \\
\tilde G^{-1} K & = \left( P_x \otimes P_y \right) \cdot \left( E_x D_x \otimes E_y + E_x \otimes E_y D_y \right) \cdot \left( P_x^{-1} \otimes P_y^{-1} \right). \\
\end{aligned}
\end{equation}
where $G=M+{\eta}K$ with $\eta=\frac{\tau^2}{4}$, and
\begin{equation}
E_\xi = ( I + \eta D_\xi )^{-1}, \qquad \xi = x, y.
\end{equation}

If we use the following identity:
\begin{equation}
\begin{aligned}
I & = P_x I_x P_x^{-1} \otimes P_y I_y P_y^{-1}  = \left( P_x \otimes P_y \right) \cdot \left( I_x \otimes I_y \right) \cdot \left( P_x^{-1} \otimes P_y^{-1} \right), \\
\end{aligned}
\end{equation}
then, the blocks of the amplification matrix are
\begin{equation}
\begin{aligned}
\Xi_{11} & = \left( P_x \otimes P_y \right) \cdot \left( I_x\otimes I_y - \tau^2 \left(E_x D_x \otimes E_y + E_x \otimes E_y D_y\right) \right) \cdot \left( P_x^{-1} \otimes P_y^{-1} \right), \\
\Xi_{12} &  = \left( P_x \otimes P_y \right) \cdot \left(I_x\otimes I_y-\tau^2\left(E_x D_x \otimes E_y + E_x \otimes E_y D_y\right) \right) \cdot \left( P_x^{-1} \otimes P_y^{-1} \right), \\
\Xi_{13} &  = \left( P_x \otimes P_y \right) \cdot \Big( -\frac{1}{2}(I_x\otimes I_y)+\big((I_x\otimes I_y)-E-\\&\frac{\tau^2}{2}\left(E_x D_x \otimes E_y + E_x \otimes E_y D_y\right) \big) \Big) \cdot \left( P_x^{-1} \otimes P_y^{-1} \right), \\
\Xi_{21} &  = \left( P_x \otimes P_y \right) \cdot \left( -\tau^2 \frac{1}{2}\left(E_x D_x \otimes E_y + E_x \otimes E_y D_y\right) \right) \cdot \left( P_x^{-1} \otimes P_y^{-1} \right), \\
\Xi_{22} &  = \left( P_x \otimes P_y \right) \cdot \left(  I_x\otimes I_y-\tau^2\frac{1}{2}\left(E_x D_x \otimes E_y + E_x \otimes E_y D_y\right) \right) \cdot \left( P_x^{-1} \otimes P_y^{-1} \right), \\
\Xi_{23} &  = \left( P_x \otimes P_y \right) \cdot \left(  \left(I_x\otimes I_y-\frac{1}{2}\left(E_x \otimes E_y\right)-\frac{\tau^2 }{4}\left(E_x D_x \otimes E_y + E_x \otimes E_y D_y\right)\right) \right) \cdot \left( P_x^{-1} \otimes P_y^{-1} \right), \\
\Xi_{31} &  = \left( P_x \otimes P_y \right) \cdot \left( -{\tau^2} \left(E_x D_x \otimes E_y + E_x \otimes E_y D_y\right) \right) \cdot \left( P_x^{-1} \otimes P_y^{-1} \right), \\
\Xi_{32} &  = \left( P_x \otimes P_y \right) \cdot \left( -\tau^2\left(E_x D_x \otimes E_y + E_x \otimes E_y D_y\right) \right) \cdot \left( P_x^{-1} \otimes P_y^{-1} \right), \\
\Xi_{33} &  = \left( P_x \otimes P_y \right) \cdot \left( I_x\otimes I_y-\left(E_x \otimes E_y\right)-\frac{\tau^2}{2}\left(E_x D_x \otimes E_y + E_x \otimes E_y D_y\right) \right) \cdot \left( P_x^{-1} \otimes P_y^{-1} \right).\\
\end{aligned}
\end{equation}
By denoting $\zeta= E_x D_x \otimes E_y + E_x \otimes E_y D_y $ and $\tilde{E}=E_x \otimes E_y $, we write the matrix as:
\begin{equation}\label{eq:stability}
\begin{aligned}
\Xi & =
\begin{bmatrix}
P_x\otimes P_y & \bfs{0}&\bfs{0} \\
\bfs{0} & P_x\otimes P_y  &\bfs{0}\\
\bfs{0}&\bfs{0}&P_x\otimes P_y \\
\end{bmatrix}
\begin{bmatrix}
I - \tau^2 \zeta & I-\tau^2\zeta& \frac{1}{2}I+ \left(-\tilde{E}-\frac{\tau^2}{2}\zeta \right) \\[0.2cm]
- \frac{1}{2}\tau^2  \zeta &I-\frac{1}{2}\tau^2\zeta &   I-\frac{1}{2}\tilde{E}-\frac{\tau^2}{4}\zeta \\[0.2cm]
-{\tau^2} \zeta&-\tau^2\zeta&I-\tilde{E}-\frac{\tau^2}{2}\zeta\\
\end{bmatrix}^n\\
&\begin{bmatrix}
 P_x^{-1} \otimes P_y^{-1} & \bfs{0}&\bfs{0} \\
\bfs{0} & P_x^{-1} \otimes P_y^{-1}&\bfs{0} \\
\bfs{0}&\bfs{0}& P_x^{-1} \otimes P_y^{-1}\\
\end{bmatrix}.
\end{aligned}
\end{equation}
To prove the unconditional stability of the method, we calculate its spectral radius:
\begin{equation}
\tilde{\Xi}  =
\begin{bmatrix}
I - \tau^2 \zeta & I-\tau^2\zeta& \frac{1}{2}I+ \left(-\tilde{E}-\frac{\tau^2}{2}\zeta \right) \\[0.2cm]
- \frac{1}{2}\tau^2  \zeta &I-\frac{1}{2}\tau^2\zeta &   I-\frac{1}{2}\tilde{E}-\frac{\tau^2}{4}\zeta \\[0.2cm]
-{\tau^2} \zeta&-\tau^2\zeta&I-\tilde{E}-\frac{\tau^2}{2}\zeta\\
\end{bmatrix}\\
\end{equation}
First , by defining $\sigma= \tau^2 D_\xi$, we consider the two limiting cases for $\sigma$: $\sigma\rightarrow 0$ and $\sigma \rightarrow \infty$. In the limit $\sigma\to 0$, since $D_\xi$ is diagonal, $E_\xi \to I$ and consequently, we have $\tau^2 \zeta \rightarrow 0$ and $E \to I$ . Hence, $\tilde \Xi$ becomes upper triangular with the following eigenvalues:
\begin{equation}
\lambda_1=\lambda_2=1, \qquad \lambda_3=0. 
\end{equation}
Hence, due to the equal multiplicity with the dimension of the stiffness matrix in 2D, $K$, the eigenvalues are bounded by $1$, and the method is unconditionally stable. In the case of the infinite time-step size, the matrix $\tilde \Xi$ becomes:
\begin{equation}\label{eq:xi}
\begin{aligned}
\tilde \Xi& =
\begin{bmatrix}
I&I&\frac{1}{2}I\\
\bfs{0}&I&I\\
\bfs{0}&\bfs{0}&I\\
\end{bmatrix}.\\
\end{aligned}
\end{equation}
Therefore, in the limit $\sigma \to \infty$, we obtain the eigenvalues $\lambda=1$. This analysis shows that the method is stable but not A-stable. Additionally, one can show that the scheme is stable for any finite time step size. 
\begin{remark}
The study of the unconditional stability of 3D splitting follows the same logic, but it is more involved.
\end{remark}

\section{Numerical results for scalar P-wave equation}

\begin{figure}
\includegraphics[scale=1.0]{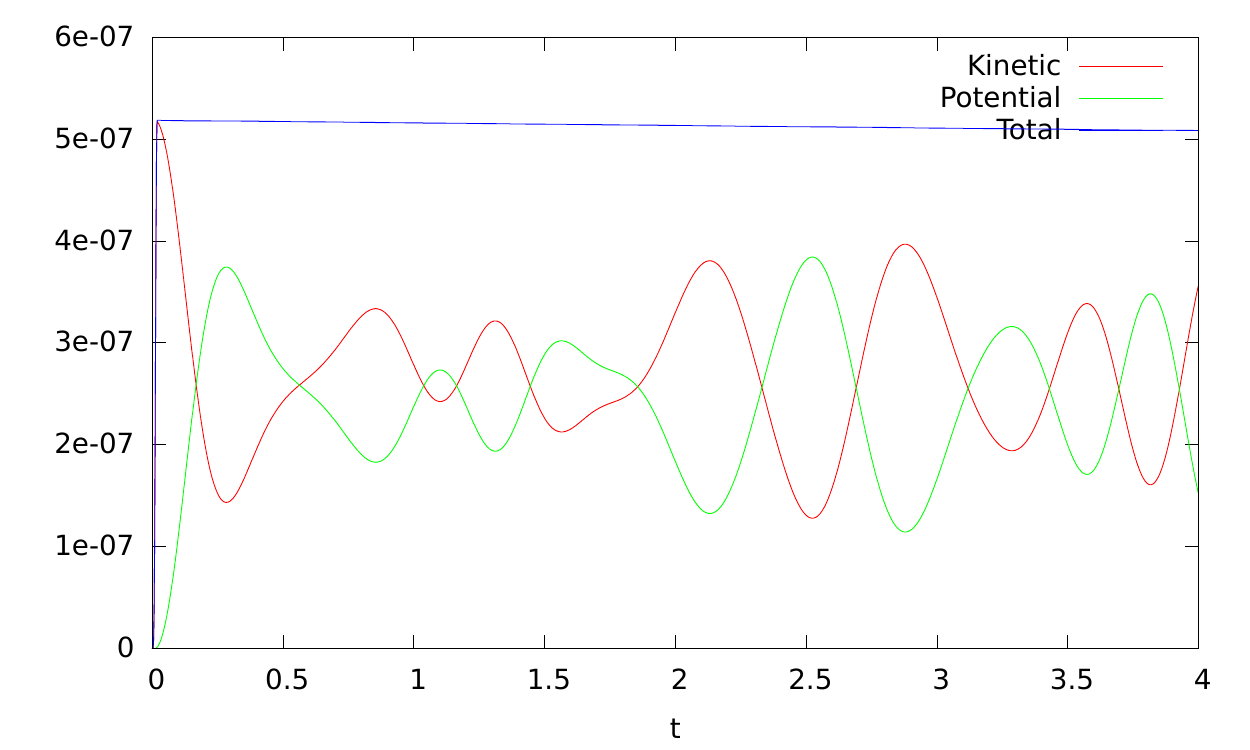}
\caption{The kinetic, potential and total energy through the entire simulation of P-wave propagation with time step $0.01$.}
\label{fig:energies}
\end{figure}

\begin{figure}
\begin{subfigure}[b]{0.2\textwidth}
\includegraphics[scale=0.18]{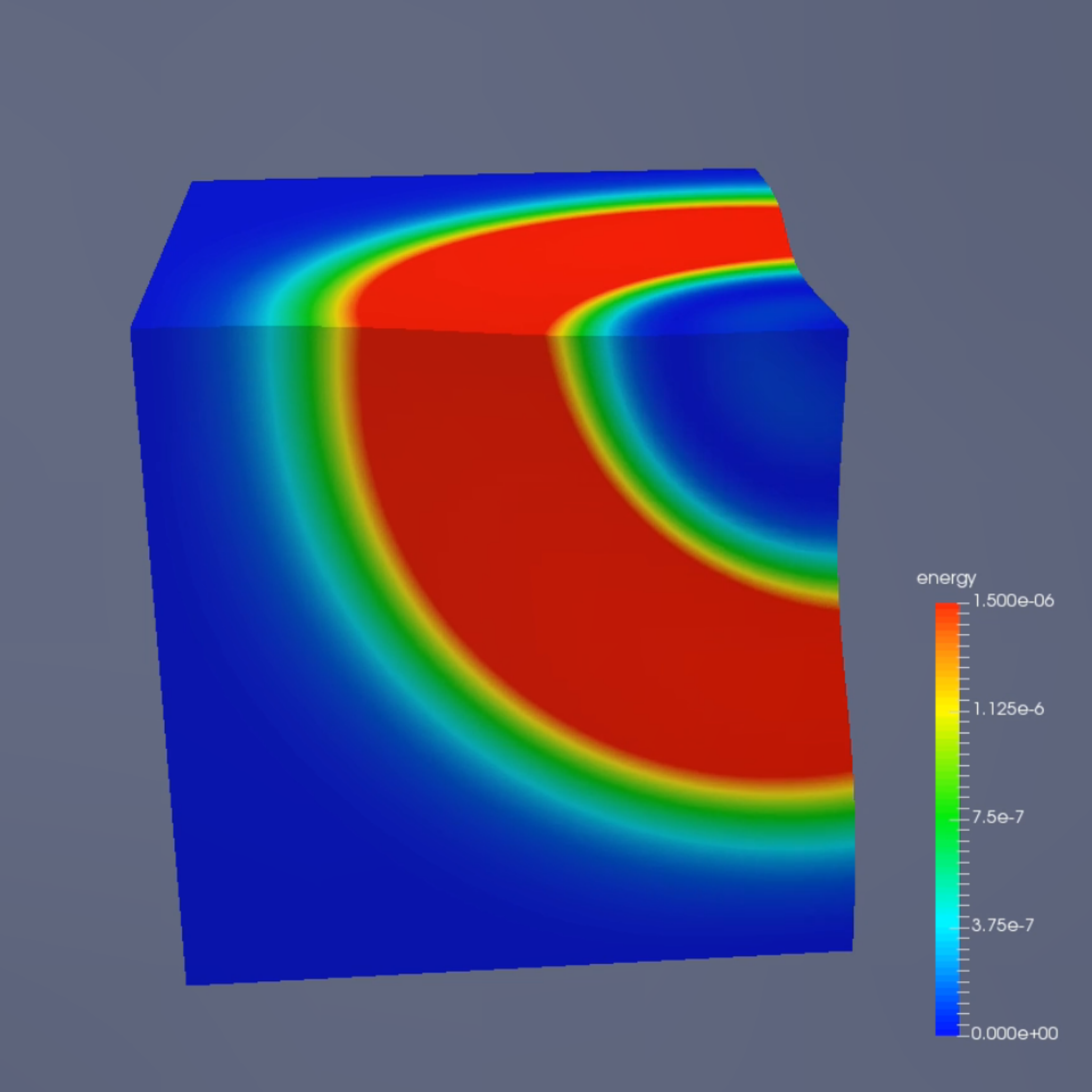}
\end{subfigure}
\begin{subfigure}[b]{0.2\textwidth}
\includegraphics[scale=0.18]{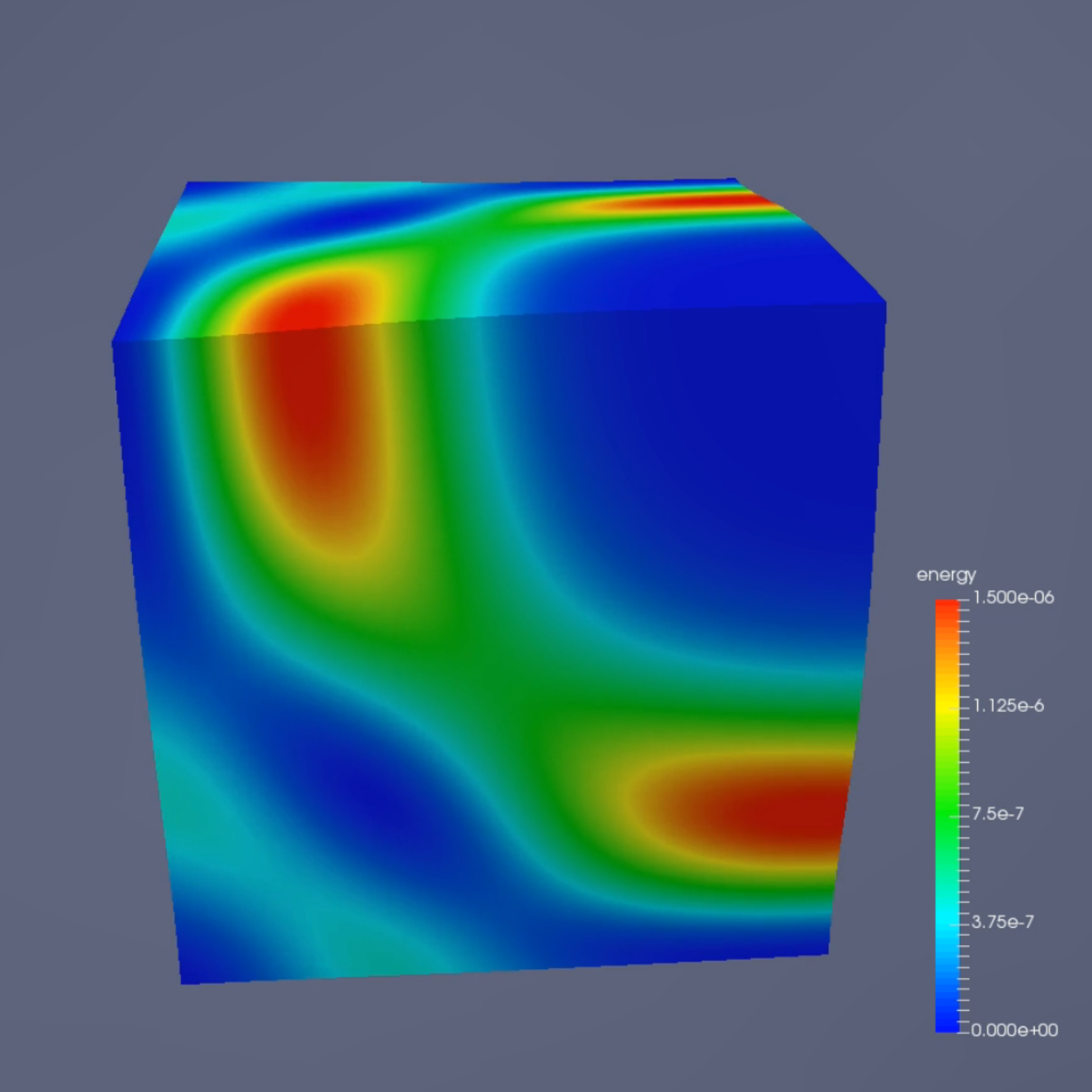}
\end{subfigure}
\begin{subfigure}[b]{0.2\textwidth}
\includegraphics[scale=0.18]{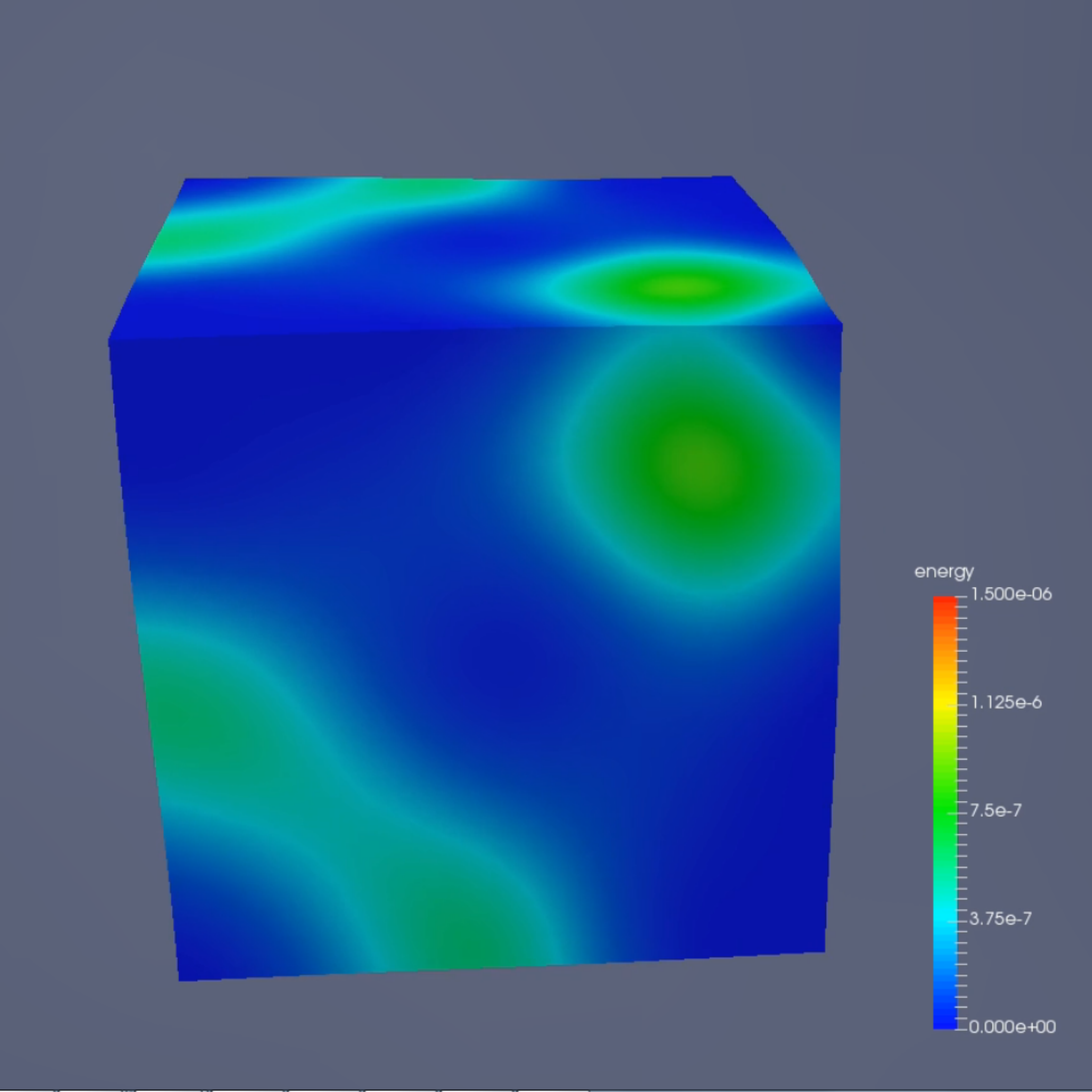}
\end{subfigure}
\begin{subfigure}[b]{0.2\textwidth}
\includegraphics[scale=0.18]{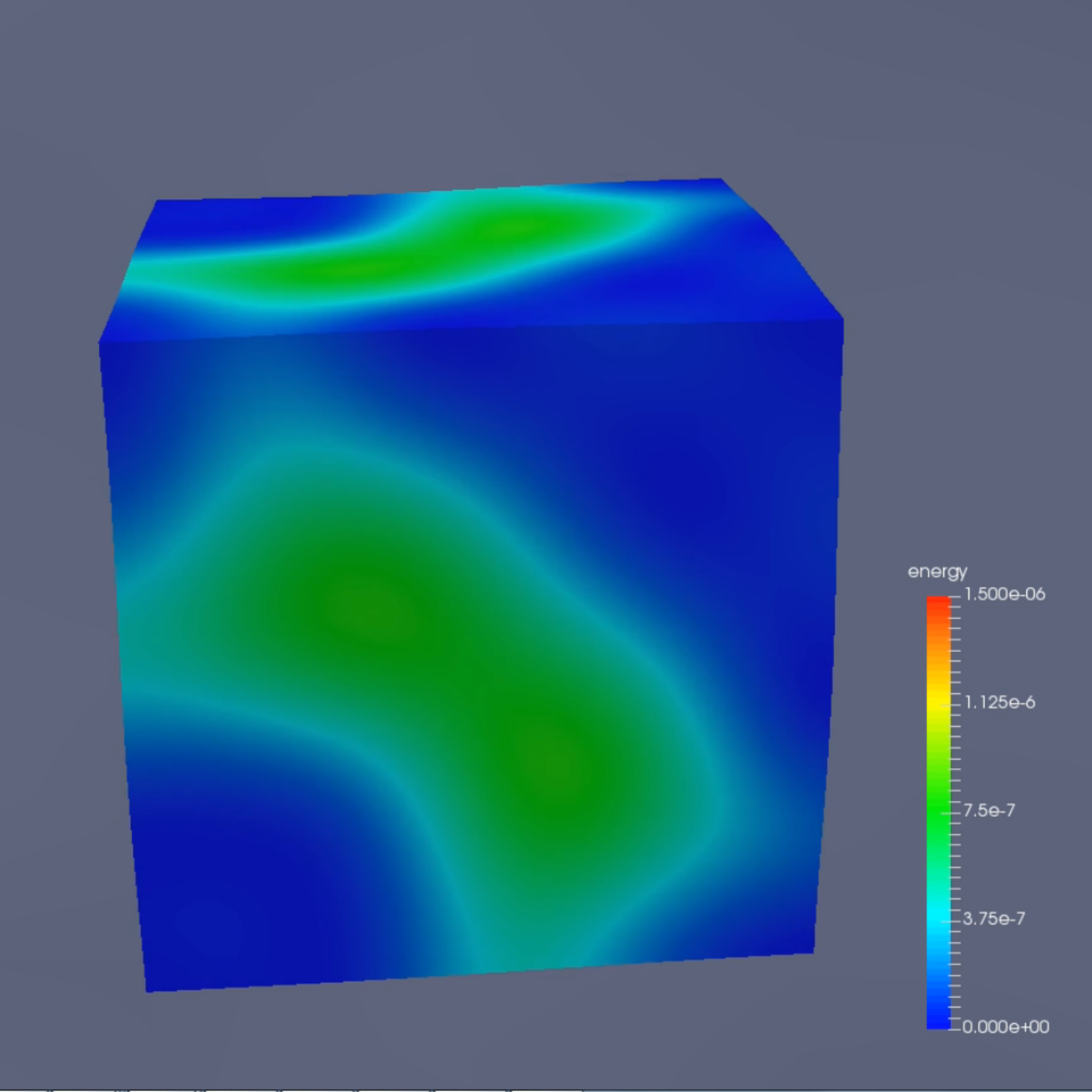}
\end{subfigure}
\caption{Stable simulation.}
\label{le_fig10e2b}
\end{figure}

We test our algorithm in a scalar P-wave propagation problem over a three-dimensional mesh with $32\times 32\times 32$ elements and time step size $dt=0.01$. We plot in Figures~\ref{fig:energies} and~\ref{le_fig10e2b} the kinetic, potential, and total energies through the entire simulation, as well as some snapshots from intermediate time steps. 

\begin{figure}
\includegraphics[scale=0.5]{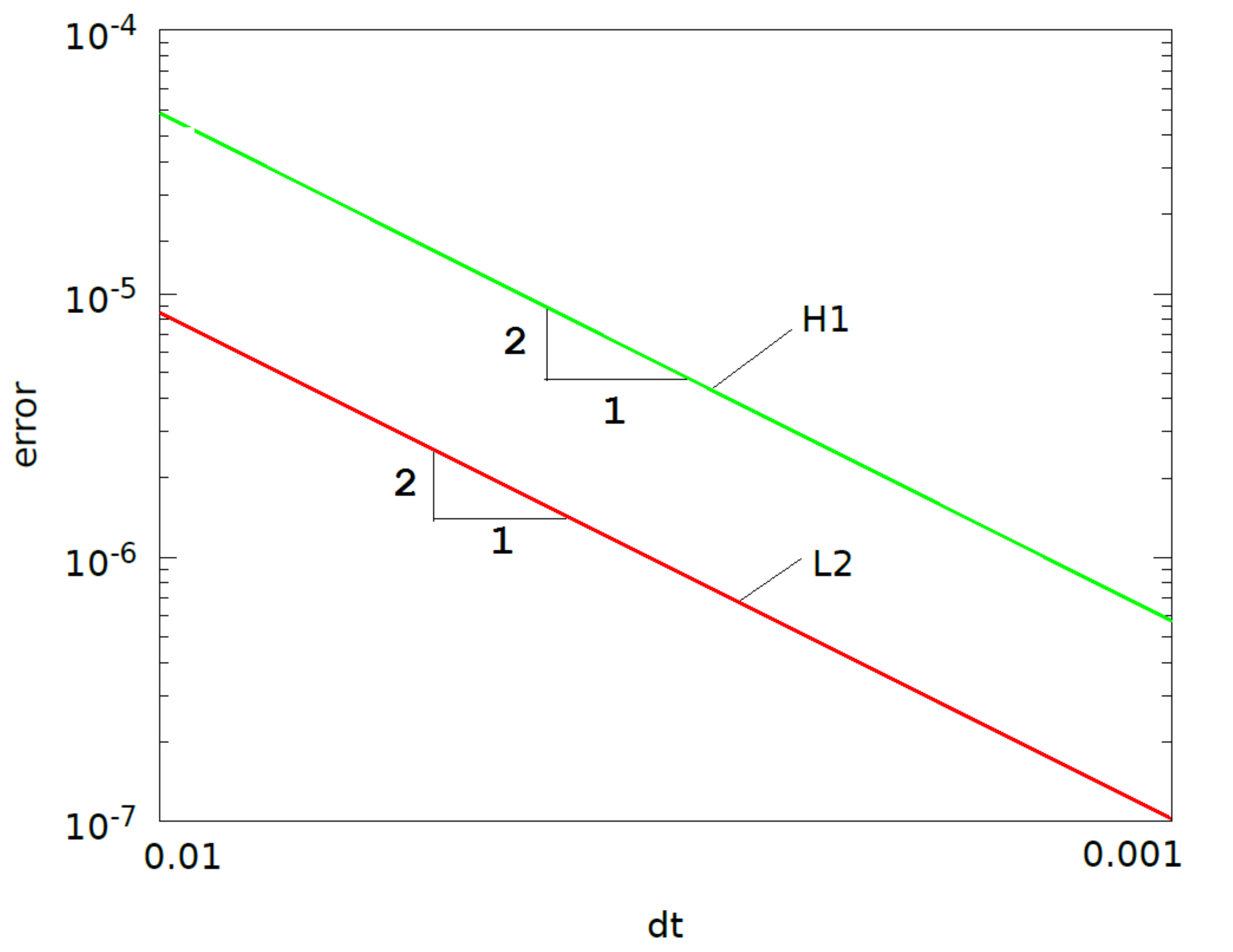}
\caption{The second order time integration scheme for P-wave equation.}
\label{fig:Pwave_order}
\end{figure}

We also verify numerically second order in time of the method, as presented in Figure~\ref{fig:Pwave_order}.

\section{Direction splitting for elastic wave propagation}

In this section, we solve the linear elasticity problem given by
\begin{equation}
  \left\{
  \begin{aligned}
    \rho\, \partial_{tt} \Vect{u} & =
      \nabla\cdot\SGM+\Vect{F} & \text{on }\Omega\times[0,T]\\
    \boldsymbol{u}(x, 0) & = u_0 & \text{for }x\in\Omega\\
    \SGM\cdot \Vect{\hat{n}} & = 0 & \mbox{on }\partial\,\Omega
  \end{aligned}
  \right.
\end{equation}
where $\Omega=[0, 1]^3$ is a unit cube, $\Vect{u}$~is a three-dimensional displacement vector to be calculated, $\rho$~is material density, $\Vect{F}$~is the applied external force, and~$\boldsymbol{\sigma}$ is the Cauchy stress tensor, given by
\begin{gather}
  \sigma_{ij}=c_{ijkl}\epsilon_{lk},
  \qquad
  \epsilon_{ij} = \frac{1}{2}\left(\partial_j u_i+ \partial_i u_j\right)
\end{gather}
and $\mathbf{c}$ is the elasticity tensor. Corresponding semi-discretized weak formulation is given by
\begin{eqnarray}
  \begin{aligned}
   \left(\ddot{u}^{n+1}_i,w\right) = \frac{1}{\rho} \left(\sigma_{ij,j}+F_i,w\right),
  \end{aligned}
\end{eqnarray}
where for repeated indexes we apply the Einstein summation convention, and 
\begin{equation}
\sigma_{ij} = 2\mu \epsilon_{ij}+\lambda \epsilon_{kk} \delta_{ij},
\end{equation}
by denoting $\epsilon_{ij} = \frac{1}{2}\left(\partial_j u_i+ \partial_i u_j\right)$. The weak form is obtained by taking the scalar product with a test functions $w_i$ and integrating by parts
\begin{equation}
\rho (w_i,\ddot{u}^{n+1}_i) + (w_{(i,j)},\sigma_{ij})=(F_i,w),
\end{equation}
where
\begin{equation}
w_{(i,j)}=\frac{w_{i,j}+w_{j,i}}{2}.
\end{equation}
We discretize
\begin{equation}
w_i = \sum_{ab} {n_am_b C^i_{ab}} \quad 
\ddot{u}^{n+1}_i = \sum_{cd} {n_cm_d D^i_{cd}} 
\end{equation}
where $n_a m_b$ and $n_c m_d$ denotes the tensor product two-dimensional B-spline basis functions, and $C^i_{ab}, i=1,2$ and $D^i_{cd}$ denotes the coefficients. We can obtain
\begin{equation}
w_{1,1} = n'_a m_b C^1_{ab}
\end{equation}
\begin{equation}
w_{1,2} = n_a m'_b C^1_{ab}
\end{equation}
\begin{equation}
w_{2,1} = n'_a m_b C^2_{ab}
\end{equation}
\begin{equation}
w_{2,2} = n_a m'_b C^2_{ab}
\end{equation}
Thus,
\begin{eqnarray}
 w_{(1,1)} = \frac{w_{1,1} + w_{1,1}}{2} =  w_{1,1}=n'_a m_b C^1_{ab}, \nonumber \\
w_{(1,2)} = \frac{w_{1,2} + w_{2,1}}{2} = \frac{n_a m'_b C^1_{ab} + n'_a m_b C^2_{ab}}{2}, \nonumber \\
 w_{(2,1)} = \frac{w_{2,1} + w_{1,2}}{2} = w_{(1,2)} = \frac{n_a m'_b C^1_{ab} + n'_a m_b C^2_{ab}}{2}, \nonumber \\
w_{(2,2)} = \frac{w_{2,2} + w_{2,2}}{2} = w_{2,2} = n_a m'_b C^2_{ab}.
\end{eqnarray}
Moreover,
\begin{eqnarray}
2w_{(1,2)}=n_a m'_b C^1_{ab} + n'_a m_b C^2_{ab}.
\end{eqnarray}
We substitute the constitutive law into the weak form
\begin{equation}
(w_{(i,j)},\sigma_{ij}) = 2\mu (w_{(i,j)},\epsilon_{ij})+\lambda (w_{(i,j)},\epsilon_{kk} \sigma_{ij}).
\end{equation}
Since $u_{(i,j)}=\epsilon_{ij}$ and $u_{k,k}=u_{(k,k)}=\epsilon_{kk}$, we utilize the definition of the Kronecker delta
\begin{equation}
(w_{(i,j)},\sigma_{ij}) = 2\mu (w_{(i,j)},u_{(i,j)})+\lambda (w_{(j,j)},u_{(k,k)}).
\end{equation}
Let us rewrite the differential operator $\Upsilon$ that corresponds the linear-elasticity problem in 2D as
\begin{equation}\label{eq:lin}
\begin{aligned}
\Upsilon& =
\begin{bmatrix}
\Upsilon_{11}&\Upsilon_{12}\\
\Upsilon_{21}&\Upsilon_{22}\\
\end{bmatrix},\\
\end{aligned}
\end{equation}
where 
\begin{equation}
	\begin{aligned}
	\Upsilon_{11} &=(2\mu+\lambda) K_x \otimes M_y+\mu K_x \otimes M_y ,\\ 
		\Upsilon_{12} &={\mu}B_x \otimes B^T_y+{\lambda}B^T_x \otimes B_y,\\ 
				\Upsilon_{21} &={\mu}B^T_x \otimes B_y+{\lambda}B_x \otimes B^T_y,\\ 
			\Upsilon_{22} &= \mu K_x \otimes M_y+(2\mu+\lambda) K_x \otimes M_y,\\ 
	\end{aligned}
\end{equation}
and we also denote the mixed matrices as
\begin{eqnarray}
(n'_a,n_c)_x=B_x, \nonumber \\
(n_a,n'_c)_x=B^T_x, \nonumber \\
(n'_b,n_d)_y=B_y, \nonumber \\
(n_b,n'_d)_y=B^T_y .
\end{eqnarray}
Next, we use the idea of alternating triangular methods~\cite{ samarskii1964economical} to the first-order evolutionary equations to construct an alternative to the second-order equations where to apply our scheme~\cite{ lisbona2001operator}. The alternating triangular method allows us to extend the operator splitting given by
\begin{equation}
	\Upsilon=\Upsilon^{(1)}+\Upsilon^{(2)},
\end{equation}
where, taking into account~\eqref{eq:lin}, we define
\begin{equation}\label{eq:decom}
\begin{aligned}
\Upsilon^{(1)} =
\begin{bmatrix}
\frac{1}{2}\Upsilon_{11}&\bfs{0}\\
\Upsilon_{21}&\frac{1}{2}\Upsilon_{22}\\
\end{bmatrix}, \qquad \qquad
\Upsilon^{(2)} =\begin{bmatrix}
\frac{1}{2}\Upsilon_{11}&\Upsilon_{12}\\
\bfs{0}&\frac{1}{2}\Upsilon_{22}\\
\end{bmatrix}.\\
\end{aligned}
\end{equation}
Finally, we solve the fully discrete problem using a two-stage approach. The predictor stage calculates $\tilde{U}^{n+1}=\left[\tilde{U}_x^{n+1},\tilde{U_y}^{n+1}\right]^T$ as 
\begin{equation}\label{eq:pre}
\rho M_x \otimes M_y\left( \frac{\tilde{U}^{n+1}-2U^n+U^{n-1}}{\tau^2}  \right)+\Upsilon^{(1)}\left( \frac{\tilde{U}^{n+1}+U^{n-1}}{2}  \right)+\Upsilon^{(2)}U^n=f^n.
\end{equation}  
To enhance the solution, we solve the following corrector stage
\begin{equation}\label{eq:corr}
\rho M_x \otimes M_y\left( \frac{{U}^{n+1}-2U^n+U^{n-1}}{\tau^2}  \right)+\Upsilon^{(1)}\left( \frac{\tilde{U}^{n+1}+U^{n-1}}{2}  \right)+\Upsilon^{(2)}\left( \frac{{U}^{n+1}+U^{n-1}}{2}  \right)=f^n.
\end{equation}  
Following the approach of~\eqref{eq:pre} and~\eqref{eq:corr}, one can solve two uncoupled problems to find $\tilde{U}_x^{n+1}$ and then $\tilde{U}_y^{n+1}$. Next, the corrected solution ${U}_y^{n+1}$ is obtained and is employed to find ${U}_x^{n+1}$.  To adapt the idea of splitting, we collect the terms and approximate the operators applied on the unknown vectors $\tilde{U}_x^{n+1}, \tilde{U}_y^{n+1}$ as
\begin{equation}\label{eq:spill}
\begin{aligned}
\left[M_x \otimes M_y+\frac{\tau^2}{4\rho} \left( (2\mu+\lambda) K_x \otimes M_y+\mu K_x \otimes M_y \right)\right]\tilde{U}_x^{n+1} &\simeq\\
 \left(M_x+\frac{\tau^2}{4\rho}  (2\mu+\lambda) K_x \right) \otimes \left(M_y+\frac{\tau^2}{4\rho} \mu K_y \right)\tilde{U}_x^{n+1} &,\\
 \left[M_x \otimes M_y+\frac{\tau^2}{4\rho} \left( \mu K_x \otimes M_y+(2\mu+\lambda) K_x \otimes M_y \right)\right]\tilde{U}_y^{n+1} &\simeq\\
 \left(M_x+\frac{\tau^2}{4\rho}  \mu K_x \right) \otimes \left(M_y+\frac{\tau^2}{4\rho}(2\mu+\lambda) K_y \right)\tilde{U}_y^{n+1} &.\\
\end{aligned}
\end{equation}
The splitting of the operators in the corrector stage follows the same argument as~\eqref{eq:spill}. Considering the splitting method, we approximate $M+\frac{\tau^2}{4}{\Upsilon}^{(1)}$ and $M+\frac{\tau^2}{4}{\Upsilon}^{(2)}$ using $M+\frac{\tau^2}{4}\tilde{\Upsilon}^{(1)}$ and $M+\frac{\tau^2}{4}\tilde{\Upsilon}^{(2)}$, respectively, by ignoring the higher order terms $\mathcal{O}(\tau^4)$.

\subsection{Stability of the method}
In this section, we study the stability of the resulting scheme~\eqref{eq:pre}-\eqref{eq:spill} by rewriting the generalized form as
\begin{equation}\label{eq:gen}
D \frac{{U}^{n+1}-2U^n+U^{n-1}}{\tau^2} +\Upsilon U^n=f^n, \qquad n=1,2, \cdots, N,
\end{equation}
Taking into account the decomposition~\eqref{eq:pre},~\eqref{eq:corr} and the splitting~\eqref{eq:spill}, we denote $D$ as $D=\left(\rho M +\sigma \tau^2 \tilde{\Upsilon}^{(1)}\right)\frac{1}{\rho M}\left(\rho M +\sigma \tau^2 \tilde{\Upsilon}^{(2)}\right)$. Hence, we have
\begin{equation}
	 \rho M+\sigma \tau^2\Upsilon \leq D.
\end{equation}
Here, we employ the argument proposed in~\cite{lisbona2001operator???} to study the stability as follows.
\begin{theorem}
For the method described in~\eqref{eq:pre}-\eqref{eq:spill}, the \textit{a priori estimate} for $\sigma \geq 0.25$ holds, 
\begin{equation}
\| U^{n+1} \|_* \leq \| U^{n} \|_* +\tau \|f^n\|,
\end{equation}
where
\begin{equation}
\left\| U^{n+1} \right\|_*^2 =\left\| \frac{U^{n+1}-U^{n}}{\tau} \right\|_D^2+\left \| \frac{U^{n+1}+U^{n}}{2} \right\|_\Upsilon^2.
\end{equation}
\begin{proof}
To prove, we closely follow the proof in~\cite{ lisbona2001operator} by doing the inner product of~\eqref{eq:gen} by  $\left({U}^{n+1}-2U^n+U^{n-1}\right)$. We obtain
\begin{equation}\label{eq:equality}
\begin{aligned}
&\left\| \frac{U^{n+1}-U^{n}}{\tau} \right\|_D^2-\left\| \frac{U^{n}-U^{n-1}}{\tau} \right\|_D^2+\left\| \frac{U^{n+1}+U^{n}}{2} \right\|_\Upsilon^2-\left\| \frac{U^{n}+U^{n-1}}{2} \right\|_\Upsilon^2\\
&= \left(f^n,{U}^{n+1}-2U^n+U^{n-1}\right)
\end{aligned}
\end{equation}
The left-hand side of the~\eqref{eq:equality} becomes
\begin{equation}
	\left(\left\| {{U}^{n+1}} \right \|_*+\left\| {{U}^{n}} \right \|_*\right)\left(\left\| {{U}^{n+1}} \right \|_*-\left\| {{U}^{n}} \right \|_*\right)
\end{equation}
On the right-hand side of~\eqref{eq:equality}, we have
\begin{equation}
\begin{aligned}
\left(f^n,\frac{{U}^{n+1}-2U^n+U^{n-1}}{\tau}\right) &\leq \left\| \frac{{U}^{n}-U^{n-1}}{\tau} \right \|\left(\left\| \frac{{U}^{n}-U^{n-1}}{\tau} \right \|+\left\| \frac{{U}^{n+1}-U^{n}}{\tau} \right \|\right)\\
&\leq \left\| f^n \right \| \left(\left\| {{U}^{n+1}} \right \|_*+\left\| {{U}^{n}} \right \|_*\right)
\end{aligned}
\end{equation}
Hence, we obtain an \textit{a priori} estimate for the method that establishes its stability with respect to the initial data and the right-hand side. This completes the proof.
\end{proof}
\end{theorem}
\begin{remark}
Again, for the sake of brevity, we omit the proof for 3D problems, which follows the same logic. 
\end{remark}

\section{Numerical results for linear elasticity}

\begin{figure}
\includegraphics[scale=1.0]{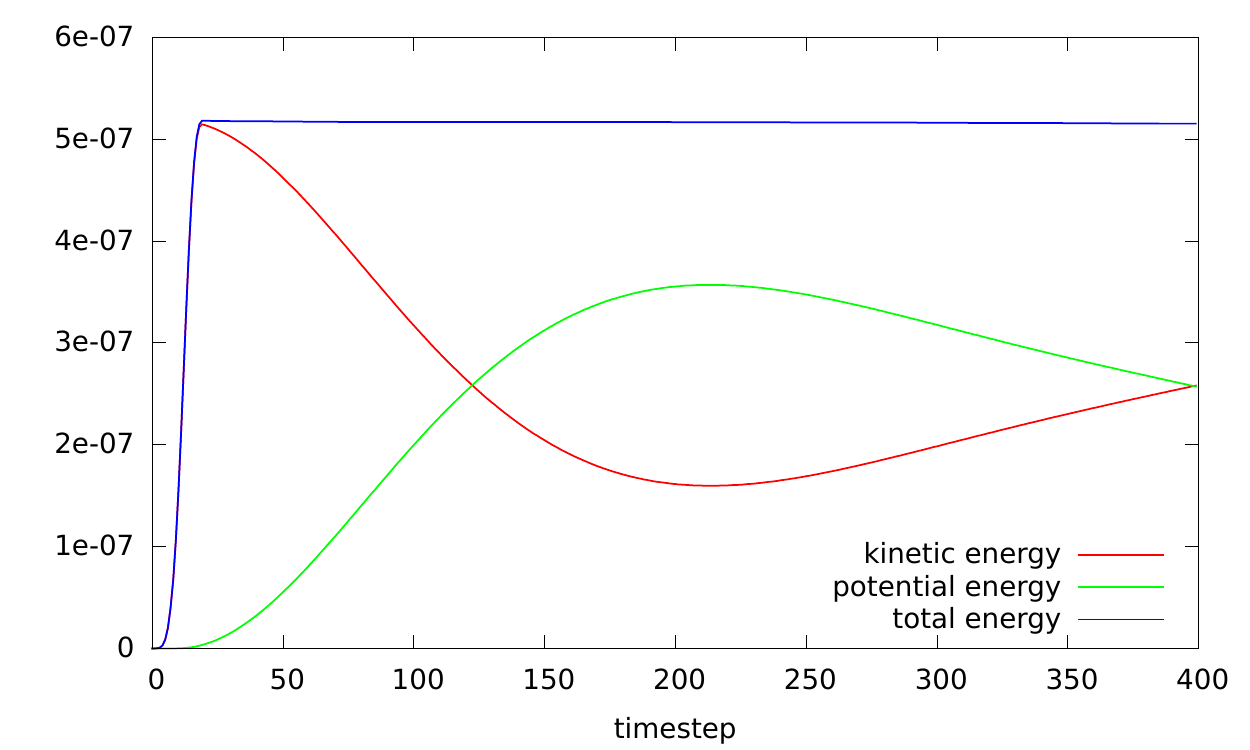}
\caption{The kinetic, potential, and total energies for the entire simulation, time-step size $10^{-2}$. The total energy remains constant.}
\label{fig:energiesLE}
\end{figure}

\begin{figure}
\includegraphics[scale=0.5]{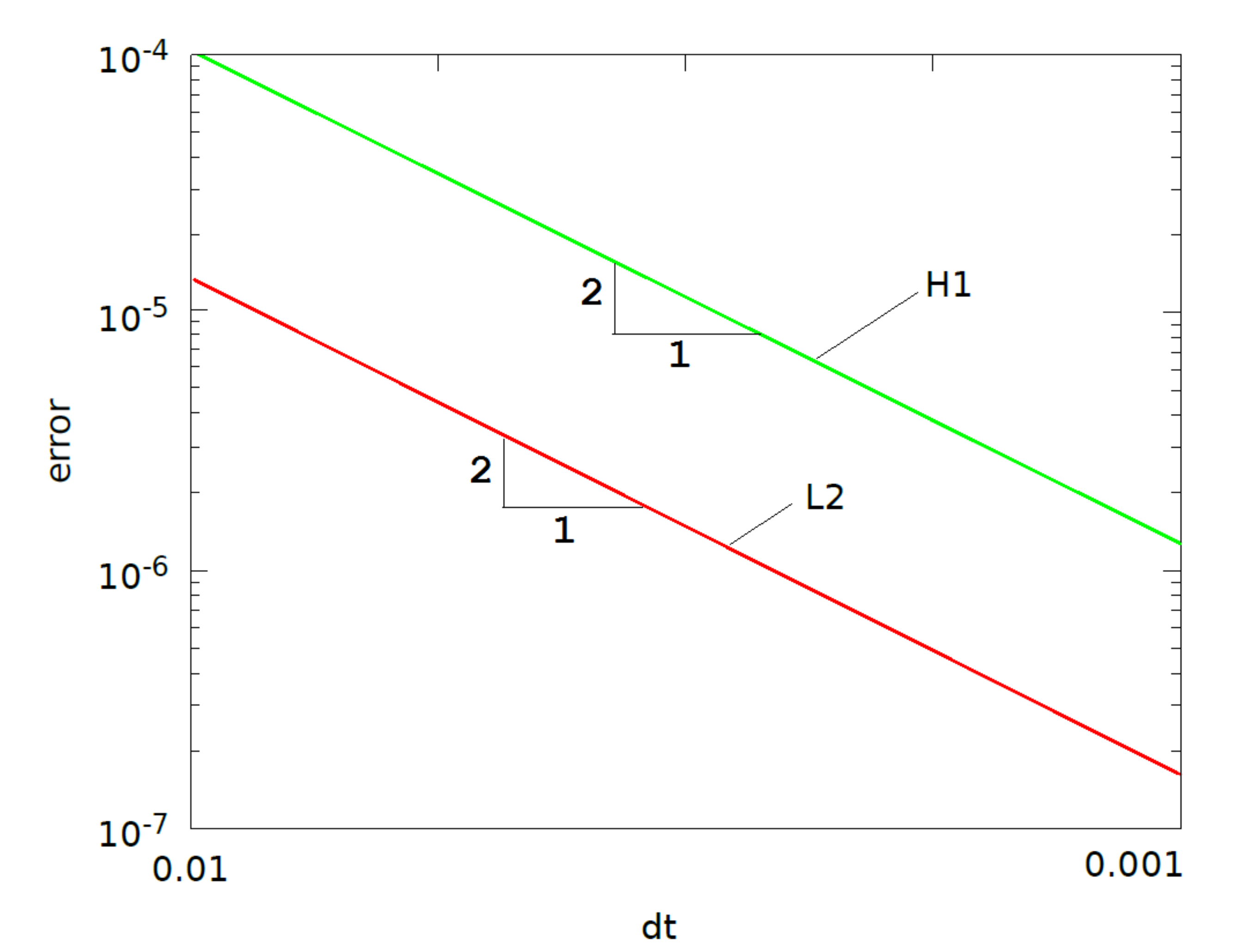}
\caption{The second order time integration scheme for linear elasticity.}
\label{fig:LE_order}
\end{figure}

\begin{figure}
\begin{subfigure}[b]{0.2\textwidth}
\includegraphics[scale=0.185]{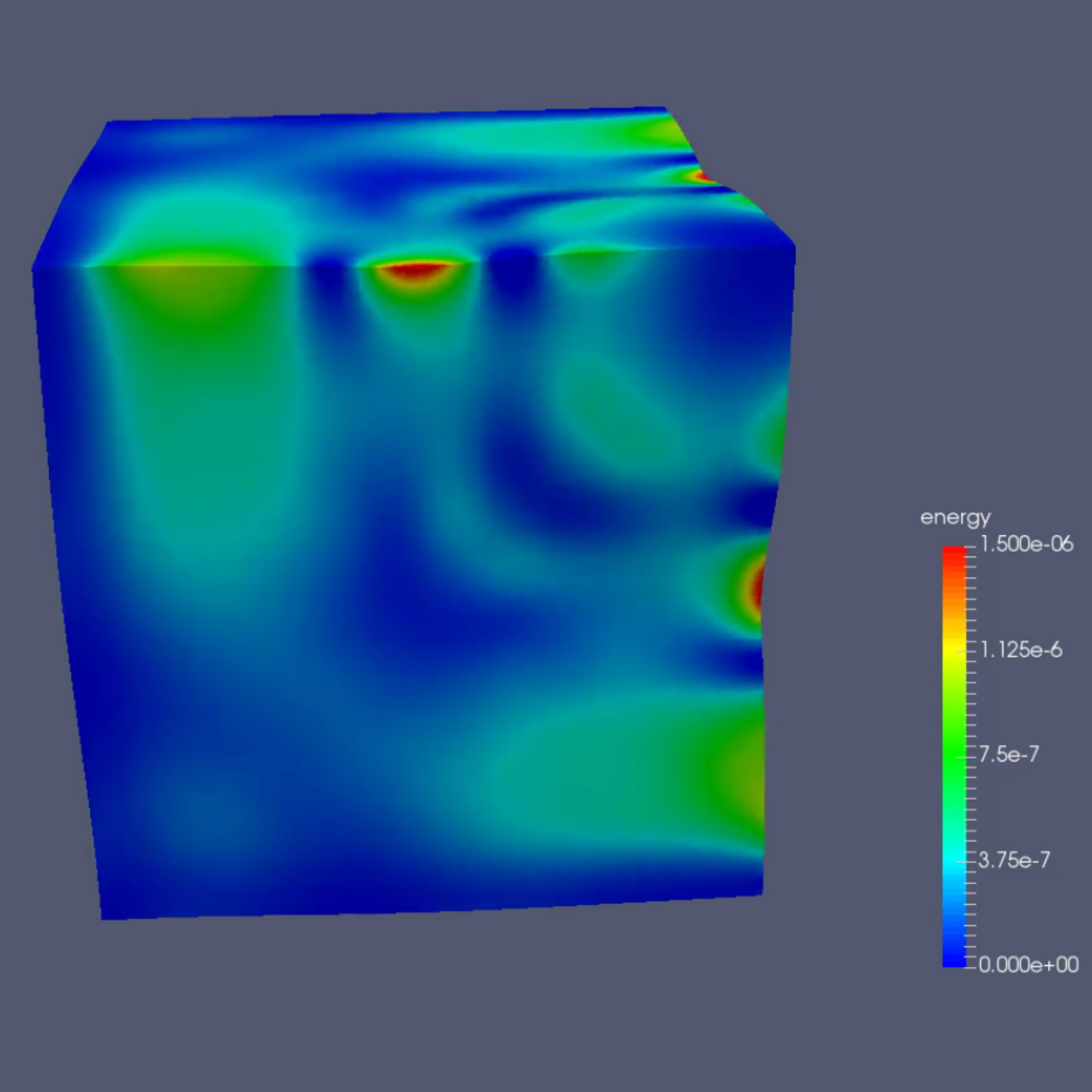}
\end{subfigure}
\begin{subfigure}[b]{0.2\textwidth}
\includegraphics[scale=0.185]{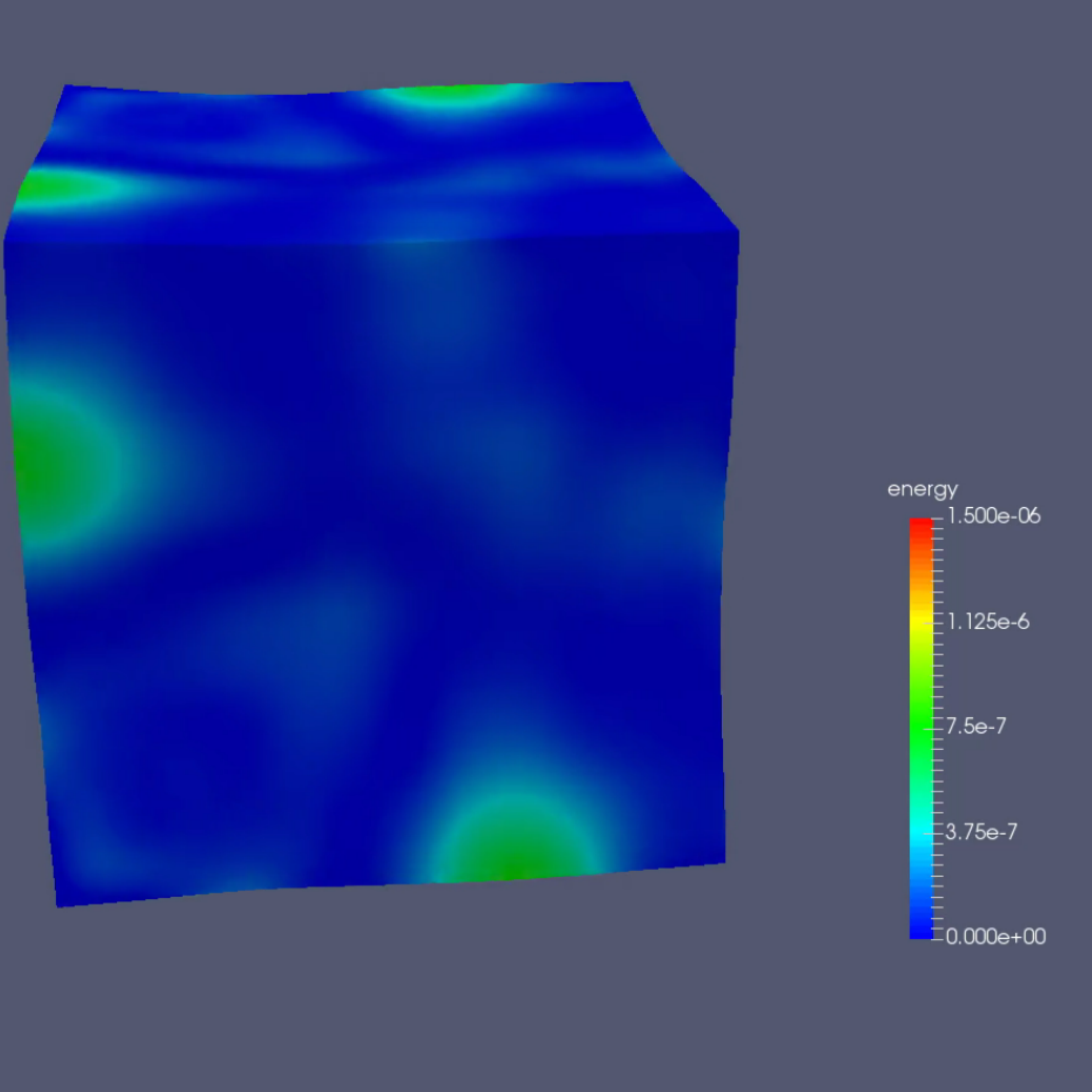}
\end{subfigure}
\begin{subfigure}[b]{0.2\textwidth}
\includegraphics[scale=0.185]{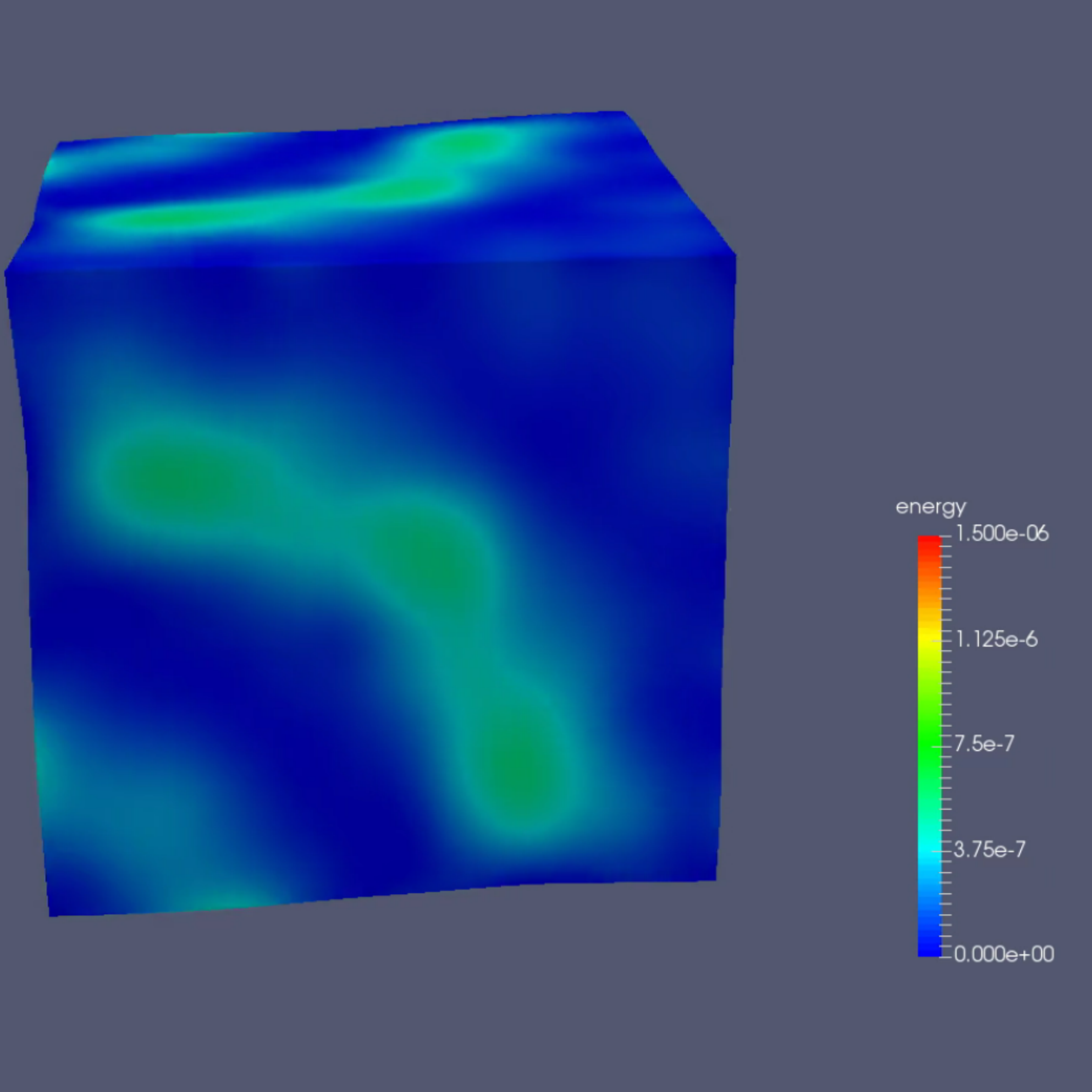}
\end{subfigure}
\begin{subfigure}[b]{0.2\textwidth}
\includegraphics[scale=0.185]{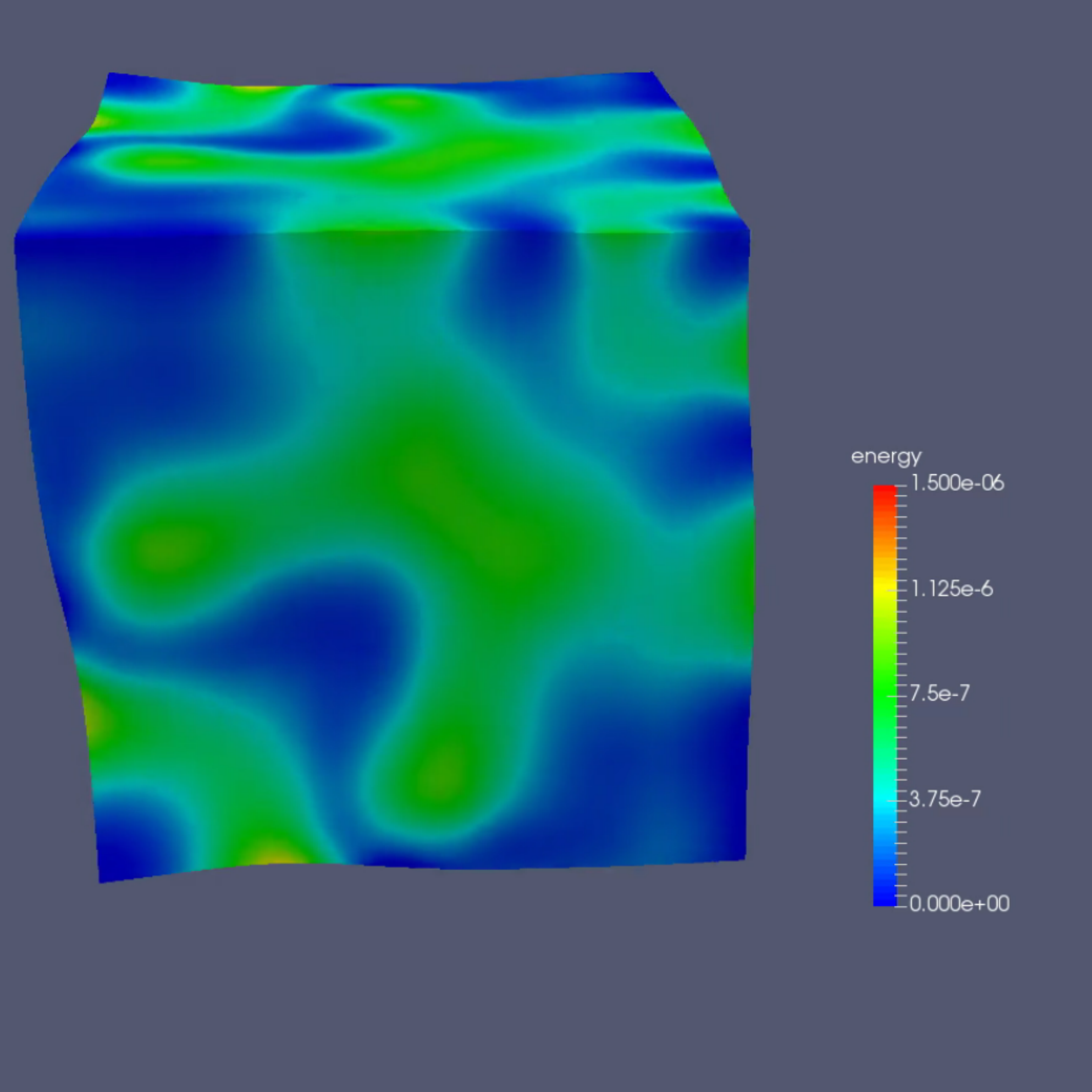}
\end{subfigure}
\caption{Stable simulation.}
\label{le_fig10e2bLE}
\end{figure}

We apply our algorithm to a linear elasticity problem in 3D with a mesh composed of $32^3$ elements with a time-step size $10^{-2}$. We plot the evolution of the kinetic, potential, and total energies through the entire simulation. We also provide snapshots from intermediate time steps, see Figures~\ref{fig:energiesLE} and~\ref{le_fig10e2bLE}. We  verify numerically that the method has second-order accuracy in time of the method, see Figure~\ref{fig:LE_order}.

\section{Conclusions}
\label{sec:conc}

In this paper, we introduce a space-time discretization using the alternating direction method to simulate hyperbolic problems.  In particular, we use high-order, smooth isogeometric basis functions in space and an implicit time marching scheme in time. We build the spatial discretization on tensor-product spaces. We then use the Kronecker-product structure of the algebraic system to invert a sequence of implicit time steps with a cost proportional to the total number of degrees of freedom in the system. We analyze the stability of the hyperbolic solvers theoretically and, then, verify the results numerically. Namely, we show the unconditional stability of the methods and verify the second-order accuracy of the time scheme experimentally. We show the performance for 2D and 3D for the scalar and vectorial differential systems. Future work will involve development of splitting schemes for Maxwell equations \cite{Maxwell} and performing parallel version of the code \cite{Cost}.

\subsection*{Acknowledgments}

National Science Centre, Poland, partially funded the work of Maciej Paszy\'{n}ski, Marcin \L{}o\'{s}, and the visit of Pouria Behnoudfar to Krak\'{o}w via the grants 2017/26/M/ST1/00281 and 2015/ 19/B/ST8/01064. This publication was also made possible in part by the CSIRO Professorial Chair in Computational Geoscience at Curtin University and the Deep Earth Imaging Enterprise Future Science Platforms of the Commonwealth Scientific Industrial Research Organisation, CSIRO, of Australia. The European Union's Horizon 2020 Research and Innovation Program of the Marie Sk\l{}odowska-Curie grant agreement No. 777778 provided additional support. At Curtin University, The Institute for Geoscience Research (TIGeR) and by the Curtin Institute for Computation, kindly provide continuing support.

\section{Appendix: Linear computational cost solver}

The matrix $\mathcal{M} = \mathcal{M}^x \otimes \mathcal{M}^y$ has a Kronecker-product structure. Each of the matrices $\mathcal{M}^\xi$ corresponds to the one-dimensional mass matrix in the direction $\xi$. In this case, we can factorize the problem with linear cost with respect to the total number of degrees of freedom in the system.

These one-dimensional matrices have entries that correspond to the integrals of the multiplication of the one-dimensional B-spline basis functions. These B-spline basis functions have local support over $p+1$ elements, so the one-dimensional matrices $\mathcal{M}^x$, $\mathcal{M}^y$ have a banded structure.
\begin{equation}
\mathcal{M}^x_{ij} = 0 \iff |i - j| > p 
\end{equation}
\begin{equation*}
	\begin{bmatrix}
    \mathcal{M}^x_{11}  & \mathcal{M}^x_{12} & \mathcal{M}^x_{13} & \mathcal{M}^x_{14} & 0 & 0 & \cdots & 0 \\
    \mathcal{M}^x_{21} & \mathcal{M}^x_{22} & \mathcal{M}^x_{23} & \mathcal{M}^x_{24} & \mathcal{M}^x_{25} & 0 & \cdots & 0 \\
    \mathcal{M}^x_{31} & \mathcal{M}^x_{32} & \mathcal{M}^x_{33} & \mathcal{M}^x_{34} & \mathcal{M}^x_{35} & \mathcal{M}^x_{36} & \cdots & 0 \\
    \vdots & \vdots & \vdots & \vdots &  \vdots & \vdots &  & \vdots\\
    0 & 0 & \ldots & \ldots & \mathcal{M}^x_{n(n-3)}& \mathcal{M}^x_{n(n-2)} & \mathcal{M}^x_{n(n-1)} & \mathcal{M}^x_{nn}
  \end{bmatrix}
\end{equation*}where $\mathcal{M}^x_{ij} = \Prod{B^x_i}{B^x_j}$. Same applies for $\mathcal{M}^y_{ij}$.

The Kronecker product structure of the matrix allows us to perform the following trick.  Rather than solving a 3D problem, we can solve three one-dimensional problems with multiple right-hand-sides.
\begin{eqnarray}
\begin{aligned}
	\begin{bmatrix}
    \mathcal{M}^x_{11}  & \mathcal{M}^x_{12} & \mathcal{M}^x_{13} & \mathcal{M}^x_{14} & 0 &  \cdots & 0 \\
    \mathcal{M}^x_{21} & \mathcal{M}^x_{22} & \mathcal{M}^x_{23} & \mathcal{M}^x_{24} & \mathcal{M}^x_{25}  & \cdots & 0 \\
    \vdots & \vdots  & \vdots &  \vdots & \vdots &  & \vdots\\
    0 & \ldots &  0 & \mathcal{M}^x_{n(n-3)}& \mathcal{M}^x_{n(n-2)} & \mathcal{M}^x_{n(n-1)} & \mathcal{M}^x_{nn}
  \end{bmatrix} \nonumber \\
    \begin{bmatrix}
    x_{111} & x_{121} & \cdots & x_{1lm}
    \\
    x_{211} & x_{221} & \cdots & x_{2lm}
    \\
    \vdots & \vdots & \ddots & \vdots \\
    x_{k11} & x_{k21} & \cdots & x_{klm}
    \\
  \end{bmatrix}
  =
  \begin{bmatrix}
    b_{111} & b_{121} & \cdots & b_{1lm}
    \\
    b_{211} & b_{221} & \cdots & b_{2lm}
    \\
    \vdots & \vdots & \ddots & \vdots \\
    b_{k11} & b_{k21} & \cdots & b_{klm}
  \end{bmatrix}
\end{aligned}
\end{eqnarray}
\begin{eqnarray}
\begin{aligned}
	\begin{bmatrix}
    \mathcal{M}^y_{11}  & \mathcal{M}^y_{12} & \mathcal{M}^y_{13} & \mathcal{M}^y_{14} & 0  & \cdots & 0 \\
    \mathcal{M}^y_{21} & \mathcal{M}^y_{22} & \mathcal{M}^y_{23} & \mathcal{M}^y_{24} & \mathcal{M}^y_{25} &  \cdots & 0 \\
    \vdots & \vdots & \vdots &  \vdots & \vdots &  & \vdots\\
   0 & \ldots & 0 & \mathcal{M}^y_{n(n-3)}& \mathcal{M}^y_{n(n-2)} & \mathcal{M}^y_{n(n-1)} & \mathcal{M}^y_{nn}
  \end{bmatrix}\nonumber \\
  \begin{bmatrix}
    y_{111} & y_{211} & \cdots & y_{k1m}
    \\
    y_{121} & y_{211} & \cdots & y_{k2m}
    \\
    \vdots & \vdots & \ddots & \vdots \\
    y_{1l1} & y_{1l1} & \cdots & y_{klm}
    \\
  \end{bmatrix}
  =
  \begin{bmatrix}
    x_{111} & x_{111} & \cdots & x_{k1m}
    \\
    x_{121} & x_{211} & \cdots & x_{k2m}
    \\
    \vdots & \vdots & \ddots & \vdots \\
    x_{1l1} & x_{2l1} & \cdots & x_{klm}
  \end{bmatrix}
\end{aligned}
\end{eqnarray}
\begin{eqnarray}
\begin{aligned}
	\begin{bmatrix}
    \mathcal{M}^z_{11}  & \mathcal{M}^z_{12} & \mathcal{M}^z_{13} & \mathcal{M}^z_{14} & 0  & \cdots & 0 \\
    \mathcal{M}^z_{21} & \mathcal{M}^z_{22} & \mathcal{M}^z_{23} & \mathcal{M}^z_{24} & \mathcal{M}^z_{25} &  \cdots & 0 \\
    \vdots & \vdots & \vdots &  \vdots & \vdots &  & \vdots\\
   0 & \ldots & 0 & \mathcal{M}^z_{n(n-3)}& \mathcal{M}^z_{n(n-2)} & \mathcal{M}^z_{n(n-1)} & \mathcal{M}^z_{nn}
  \end{bmatrix}\nonumber \\
  \begin{bmatrix}
    z_{111} & z_{121} & \cdots & z_{1l1}
    \\
    z_{212} & z_{222} & \cdots & z_{2l2}
    \\
    \vdots & \vdots & \ddots & \vdots \\
    z_{k1m} & z_{k2m} & \cdots & z_{klm}
    \\
  \end{bmatrix}
  =
  \begin{bmatrix}
    y_{111} & y_{121} & \cdots & y_{1l1}
    \\
    y_{212} & y_{222} & \cdots & y_{kl2}
    \\
    \vdots & \vdots & \ddots & \vdots \\
    y_{k1m} & y_{k2m} & \cdots & y_{klm}
  \end{bmatrix}
\end{aligned}
\end{eqnarray}
where $\mathcal{M}^x_{ij} = \Prod{B^x_i}{B^x_j}$ and $\mathcal{M}^y_{ij} = \Prod{B^y_i}{B^y_j}$ and $\mathcal{M}^z_{ij} = \Prod{B^z_i}{B^z_j}$ . The dimensions of the first problem are $n \times n $, where $n$ is the number of B-spline basis functions along $x$ axis, and we have $ml$ right-hand-sides, where $m$ is the number of B-spline basis functions along $y$ axis, and $l$ is the number of B-spline basis functions along $z$ axis. The computational complexity of factorization of such a system is $O(n*m*l)=O(N)$~\cite{ PaszynskiBook}. We have the analogous situation in the second problem, namely $ m \times m $ system with $n*l$ right-hand-sides, which results in $O(m*n*l)=O(N)$ linear computational complexity, and in the third system we have $ l \times l $ system with $n*m$ right-hand-sides, which results in $O(l*n*m)=O(N)$ linear computational complexity.

This strategy delivers a solution to the isogeometric L2 projection problem with linear $O(N)$ computational cost.  This solution method improves on the standard direct solver cost estimates for and $O(N^2)$ in three-dimensions, see~\cite{ Collier}) for the factorization of the global problem.

\end{document}